\documentclass{amsart}

\newtheorem{theorem}{Theorem}[section]

\newtheorem{lemma}[theorem]{Lemma}

\newcommand{\ZZ}{{\mathbb Z}}
\newcommand{\RR}{{\mathbb R}}
\newcommand{\CC}{{\mathbb C}}
\newcommand{\TT}{{\mathbb T}}

\newcommand{\DD}{{\mathbb D}}

\title[ The ``action" variable is not an 
invariant
for uniqueness]{
 The ``action" variable is not an invariant
for the uniqueness in the inverse scattering
problem}

\author{A. Kheifets and P. Yuditskii}

\address{
Department of Mathematics,
The College Of
 William and  Mary,
P. O. Box 8795,  Williamsburg, Virginia
23187-8795}
 \email{ sykhei@wm.edu}

\address{Department of Mathematics,
Michigan State University,
East Lansing, MI 48824}
\email{yuditski@math.msu.edu}

\begin{document}


\date{January 18, 2002}


\begin{abstract}
We give a simple example of 
non-Ñuniqueness in the inverse scattering for
Jacobi matrices: roughly speaking $S$-matrix is
analytic. Then, multiplying a reflection
coefficient by an inner function, we repair
this matrix in such a way that it does
uniquely
determine  a Jacobi matrix of Szeg\"o
class; on the other hand the transmission
coefficient remains the same. This implies the
statement given in the title.  
 \end{abstract}

\maketitle

\section{Jacobi matrices of Szeg\"o class.
Direct scattering --- Bernstein--Szeg\"o type
theorem}

As it is well known in the theory of 
completely
integrable systems the absolute values of a
reflection coefficient have played the
role of the ``action" variables and the 
arguments
of this function have meaning of the ``angle"
variables (in this case we think on the Toda
lattice as on an integrable
system) \cite{T}. Combining results of our
previous works \cite{Kh1} (see also \cite{Kh2})
and \cite{VYu}, we give a wide set of examples
where two reflection 
coefficients, having the same
absolute values, possess completely
different properties: the first one
uniquely
determines a Jacobi matrix of Szeq\"o class
 and the second one does not. 
Note that 
the proof of the main theorem
in \cite{Kh1} (and, therefore, our
result) essentially uses  the analysis of 
\cite{K} on regularization of so called
Arov--singular 
matrix functions.

Let $J$ be a Jacobi matrix defining a bounded
self--adjoint operator on $l^2(\ZZ)$:
\begin{equation}\label{i1}
J e_n=p_n e_{n-1}+q_n e_{n}+p_{n+1} e_{n+1},
\quad n\in\ZZ,
\end{equation}
where $\{e_n\}$ is the standard basis in
$l^2(\ZZ)$, $p_n>0$. The resolvent
matrix--function is defined by the relation
\begin{equation}\label{i2}
R(z)= R(z, J)=\mathcal E^*(J-z)^{-1}
\mathcal E,
\end{equation}
 where $\mathcal E:\CC^2 \to l^2(\ZZ)$ is such
that
$$
\mathcal E\bmatrix c_{-1} \\ c_0\endbmatrix=
e_{-1} c_{-1}+e_{0} c_{0}.
$$
This matrix--function possesses an integral
representation
\begin{equation}\label{i3}
R(z)=\int\frac{d\sigma}{x-z}
\end{equation}
with a $2\times 2$ matrix--measure having
a compact support on $\RR$. $J$ is unitary
equivalent to the multiplication
operator  by an
independent variable on
$$
L_{d\sigma}^2=
\left\{f=\bmatrix f_{-1}(x) \\
f_0(x)\endbmatrix:\ \int f^*\,d\sigma\,
f<\infty\right\}.
$$
The spectrum of $J$ is called absolutely
continuous if the measure $d\sigma$ is
absolutely
continuous with respect to the Lebesgue measure
on the real axis,
\begin{equation}\label{i4}
d\sigma(x)=\rho(x)\,dx.
\end{equation}

It is natural to ask how properties of
coefficients of $J$ are reflected on its
spectral properties. One is especially
interested in $J$'s ``close" to the ``free"
matrix
$J_0$ 
with constant coefficients, $p_n=1,\ q_n=0$ (so
called Chebyshev matrix).

Let us mention that $J_0$  has the following
functional representation, besides the general
one mentioned above. The resolvent
set of $J_0$ is the domain $\bar\CC\setminus
[-2,2]$. Let $z(\zeta):\DD\to\bar\CC
\setminus [-2,2]$ be a uniformization of this
domain,
$z(\zeta)=1/\zeta+\zeta$.
With respect to the standard basis
$\{t^n\}_{n\in\ZZ}$
 in
$$
L^2=\{f(t):\ \int_{\TT}|f|^2\,dm\},
$$
the matrix of the operator of multiplication by
$z(t),\ t\in\TT$,
is the
Jacobi matrix $J_0$, since $z(t)
t^n=t^{n-1}+t^{n+1}$.

We say that $J$ with absolutely continuous
spectrum $[-2,2]$ is of Szeg\"o class if its
spectral density \eqref{i4} satisfies
\begin{equation}\label{i5}
\log\det\rho(z(t))\in L^1.
\end{equation}

\begin{theorem}\label{th1} Let $J$ be of
Szeg\"o class.
 Then
\begin{equation}\label{i6}
p_n\to 1,\ q_n\to 0, \quad n\to\pm\infty.
\end{equation}
Moreover, there exist  generalized
eigenvectors
\begin{equation}\label{i7}
\begin{split}
&p_n e^+(n-1,t)+q_n e^+(n,t)
+p_{n+1} e^+(n+1,t)= z(t) e^+(n,t)\\
&p_n e^-(-n,t)+q_n e^-(-n-1,t)
+p_{n+1} e^-(-n-2,t)= z(t) e^-(-n-1,t)
\end{split}
\end{equation}
such that the following asymptotics hold true
\begin{equation}\label{i8}
\begin{split}
s(t)e^{\pm}(n,t)=&s(t) t^n+o(1),\quad
n\to +\infty\\
s(t)e^{\pm}(n,t)=&t^n+s_\mp(t) t^{-n-1}+
o(1),\quad n\to -\infty
\end{split}
\end{equation}
 in $L^2$.
\end{theorem}

Thus the eigenvectors of $J$ behave
asymptotically as the eigenvectors of $J_0$ 
(later we make more precise statement).

The matrix formed by the coefficients of 
\eqref{i8}   
\begin{equation}\label{i9}
S(t)=\bmatrix s_-& s\\
         s&s_+
\endbmatrix(t),\qquad t\in \TT,
\end{equation}
is called the scattering matrix
of $J$. It is
unitary--valued, possesses the symmetry
property $S^*(\bar t)=S(t)$ and the following
analytic property: $s(t)$ is  
an outer function.  

In what follows every matrix--function
of the form \eqref{i9} with the above listed 
properties is called a scattering matrix.
Of course we have a good reason for this,
since {\it with every matrix $S(t)$ of this
kind  one can associate a Jacobi matrix $J$
whose scattering matrix 
(associated to $J$ according to 
Theorem \ref{th1}) is the initial
matrix--function S(t)}.
However, $S(t)$, generally speaking, does not
determine $J$ uniquely. 

To clarify all above statements we need some
notation and definitions.
First of all for
a given functions $s_\pm$ we define
the metric
\begin{equation*}
\begin{split}
||f||^2_{s_\pm}=&
\frac 1 2\left\langle
\bmatrix 1 &\overline{s_\pm(t)}\\
     s_\pm(t) & 1
\endbmatrix
\bmatrix f(t)\\ \bar t f(\bar t)
\endbmatrix,
\bmatrix f(t)\\ \bar t f(\bar t)
\endbmatrix
\right\rangle\\
=&\langle f(t)+\bar t (s_\pm f)(\bar t),
f(t)\rangle,\quad f\in L^2,
\end{split}
\end{equation*}
 and we denote by $L^2_{dm,s_\pm}$ or
$L^2_{s_\pm}$ (for shortness) the closure
of $L^2$ with respect to this new metric.

The following relations set a unitary map
from $L^2_{s_+}$ to $L^2_{s_-}$:
\begin{equation}\label{i10}
\begin{split}
\bmatrix s f^+\\ s f^-\endbmatrix
(t)=&
\bmatrix s&0\\s_+&1\endbmatrix(t)
\bmatrix f^+(t)\\ \bar t f^+(\bar t)
\endbmatrix\\
=&
\bmatrix 1&s_-\\0&s\endbmatrix(t)
\bmatrix  \bar t f^-(\bar t)\\f^-(t)
\endbmatrix.\\
\end{split}
\end{equation}
Moreover, in this case,
$$
||f^+||_{s_+}^2=||f^-||_{s_-}^2=
\frac 1 2
\{||s f^+||^2+||s f^-||^2\}.
$$
It is worth to give a scalar
variant of relations between 
$f^+\in L^2_{s_+}$ and  
${f^-\in L^2_{s_-}}$:
\begin{equation*}
 s(t)f^\mp(t)=\bar t f^\pm(\bar t) +s_\pm(t)
f^\pm(t). 
\end{equation*}

\begin{theorem} \label{th2} $J$ is a Jacobi
matrix
 of Szeg\"o class
with the spectrum $E=[-2,2]$ if and only if
$J$ possesses the scattering representation,
i.e.:
there exists a unique
 matrix--function $S(t)$ of
the form \eqref{i9} (with the 
listed properties)
and a unique pair of Fourier
transforms
\begin{equation}
\mathcal F^\pm:l^2(\ZZ)\to L^2_{s_\pm},
 \quad
(\mathcal F^{\pm} Jf)(t)=z(t)(\mathcal F^\pm
f)(t),
 \end{equation}
determining each other by the relations
\begin{equation}\label{i12}
s(t)(\mathcal F^\pm f)(t)=
\bar t (\mathcal F^\mp f)(\bar t)+
s_\mp(t)(\mathcal F^\mp f)(t),
 \end{equation}
and having the following analytic properties
\begin{equation}\label{i13}
s\mathcal F^\pm(l^2(\ZZ_\pm))\subset H^2,
\end{equation}
and asymptotic properties
\begin{equation}\label{i14}
e^\pm(n,t)=t^n+o(1)\quad\text{in}\ L^2_{s_\pm},
\quad n\to +\infty,
 \end{equation}
where
$$
e^+(n,t)=(\mathcal F^+ e_n)(t),\quad
e^-(n,t)=(\mathcal F^- e_{-n-1})(t),
$$
with $\{e_n\}$ being the standard basis
in $l^2(\ZZ)$.
\end{theorem}

We point out that asymptotic relations
\eqref{i8} and \eqref{i14} are equivalent,
moreover
\eqref{i14} directly shows that the
eigenvectors of $J_0$ asymptotically are
the eigenvectors of $J$.

\section{Uniqueness and Completeness}
Before we proceed with the uniqueness theorem,
we show how to construct at least one $J$
with the given scattering matrix $S(t)$.

Consider the space
$$
H^2_{s_+}=\text{clos}_{L^2_{s_+}} H^2,
$$
and introduce the Hankel operator $\mathcal
H_{s_+}: H^2\to H^2$,
$$
\mathcal H_{s_+} f= P_+\bar t(s_+ f)(\bar t),
\quad f\in H^2,
$$
where $P_+$ is the Riesz projection from $L^2$
onto $H^2$. This operator determines the metric
in $H^2_{s_+}$:
\begin{equation*}
\begin{split}
||f||^2_{s_+}=&\langle f(t)+\bar t(s_+ f)(\bar t),
f(t)\rangle\\
=&\langle (I+\mathcal H_{s_+}) f,
f\rangle,\quad
\forall f\in H^2.
\end{split}
\end{equation*}

\begin{theorem} \label{th3} Let $S(t)$ be a
scattering matrix, i.e., the matrix of the form
\eqref{i9} with listed properties.
Then the space $H_{s_+}^2$ is  a space of
holomorphic functions with a reproducing
kernel. Moreover,
 the reproducing vector $k_{s_+}$:
$$
\langle f, k_{s_+}\rangle=f(0),\quad\forall f\in
 H_{s_+}^2,
$$
is of the form
\begin{equation}
 k_{s_+}=(I+\mathcal H_{s_+} )^{[-1]} \text{\bf
1}:=
\lim_{\epsilon\to 0^+}(\epsilon+I+\mathcal
H_{s_+} )^{-1}
\text{\bf 1}
\quad\text{in}\ L^2_{s_+}. 
\end{equation}

Put $K_{s_+}(t)=k_{s_+}(t)/\sqrt{k_{s_+}(0)}$.
Then the system of functions
$\{t^n K_{s_+t^{2n}}(t)\}_{n\in\ZZ}$ forms an
orthonormal basis in $L^2_{s_+}$. With respect to
this basis, the  multiplication operator by
$z(t)$ is a Jacobi matrix $J=J[s_+]$ of 
Szeg\"o class.  Moreover,
the initial $S(t)$ serves as
 the scattering
matrix--function, associated with 
given $J$ by Theorem \ref{th2}, and
the relations
$$
\mathcal F^+(e_n)= t^n K_{s_+t^{2n}}(t)
$$
determine uniquely corresponding Fourier
transforms.
\end{theorem}

Let us fix the notation $J[s_+]$ for the Jacobi
matrix associated with $S(t)$ by 
Theorem \ref{th3}.
On the other hand, the system of functions
$\{t^n K_{s_-t^{2n}}(t)\}_{n\in\ZZ}$ forms an
orthonormal basis in $L^2_{s_-}$, and we can
define a Jacobi matrix 
$\tilde J=J[s_-]$ by
the relation
$$
z(t)\tilde e^+(n,t)=
\tilde p_n \tilde e^+(n-1,t)+\tilde q_n
\tilde e^+(n,t)
+\tilde p_{n+1} \tilde e^+(n+1,t),
$$
where $\{\tilde e^+(n,t)\}$ is the dual system
to the system $\{t^n K_{s_-t^{2n}}(t)\}$
(see \eqref{i12}), i.e.:
$$
s(t)\tilde e^+(-n-1,t):=\bar t^{n+1}
K_{s_-t^{2n}} (\bar t)+ s_-(t) t^n
K_{s_-t^{2n}}(t).
$$
 None guarantees 
that operators 
$J[s_+]$ and $J[s_-]$ are the same
(see beginning of the next section).
However, if
$J[s_+]=J[s_-]$, then the uniqueness
theorem  takes place.

\begin{theorem}  A scattering matrix $S(t)$
determines a
Jacobi matrix $J$ of Szeg\"o class in a 
unique way
if and only if
the following relations take place
\begin{equation}\label{i16}
s(0)K_{s_\pm}(0) K_{s_\mp t^{-2}}(0)=1.
 \end{equation}
\end{theorem}

Of course it is hard  to check identities,
especially using computer simulation, but,
in fact, \eqref{i16} has a specific
approximating meaning.

Let us return to the matrix $J[s_-]$.
According to \eqref{i13} the space
$\mathcal F^+(l^2(\ZZ_+))$ is a subspace
of $L^2_{s_+}$ consisting of holomorphic in
$\DD$ functions. In the given case we can
 even specify this space in the form
\begin{equation}\label{i17}
\mathcal F^+(l^2(\ZZ_+))=
\widehat H_{s_+}^2:=
\{f\in L^2_{s_+}:\ sf\in H^2\}.
\end{equation}
Every function $f$ from $H^2_{s_+}$ possesses
the property $sf\in H^2$, but this means
only inclusion:
\begin{equation}
 \widehat H_{s_+}^2\supset
H^2_{s_+} 
\end{equation}
The meaning of \eqref{i16} is that every
function from $\widehat H_{s_+}^2$ can be
approximated by functions from $H^2$ in
$L^2_{s_+}$--norm ($H^2$ is dense in
$\widehat H_{s_+}^2$). In
fact, it is enough  to prove that we can 
approximate
just two  functions
$\tilde e^+(0,t)$ and $\tilde e^+(1,t)$
from $\widehat H_{s_+}^2$.
This would guarantee \eqref{i16}, completeness
and the uniqueness theorem.

\section{The result}
Our first remark is
almost evident. Let us pick an analytic
scattering matrix (all entries are in
$H^\infty$). For example,
\begin{equation*}
S(t)=\bmatrix \frac{1+\Delta} 2&
\frac{1-\Delta} 2\\
         \frac{1-\Delta} 2&
\frac{1+\Delta} 2
\endbmatrix(t),
\end{equation*}
where $\Delta(t)=\overline{\Delta(\bar t)}$
is a symmetric inner function. Then the
analyticity of $s_+$ implies that 
$\mathcal H_{s_+}=0$ , 
$H^2_{s_+}=H^2$ and 
$t^n K_{s_+t^{2n}}(t)=t^n$, $n\ge 0$.
Thus, 
$$
J[s_+]|l^2(\ZZ_{+}(1))
=J_0|l^2(\ZZ_{+}(1)),
$$
where $\ZZ_{+}(m)=\{n\in\ZZ: n\ge m\}$.  
 The analyticity of $s_-$, in its turn, implies
that
$$
J[s_-]|l^2(\ZZ_{-}(-1))
=J_0|l^2(\ZZ_{-}(-1)),
$$
$\ZZ_{-}(m)=\{n\in\ZZ: n\le m-1\}$.
 
Thus, the uniqueness would imply that 
$J$ associated to given $S(t)$ coincides
with $J_0$
up to three coefficients: $p_0, q_0,
q_{-1}$. For fun we give an exact formula
for $S(t)$ related to this case:
\begin{equation}
S(t)=\frac 1{\phi(t)}
\bmatrix\psi_*(t)& p_0(1-t^2)\\
p_0(1-t^2)&\psi(t)
\endbmatrix,
\end{equation}
where $\phi(t)=(1-q_0 t)(1-q_{-1} t)-
p_0^2 t^2$,
$\psi(t)=(1-q_0 t)(q_{-1} -t)+
p_0^2 t^2$ and 
$\psi_*(t)=t^2\overline{\psi(t)}$.
 Note that this matrix has no bound states
($s(t)$ is holomorphic in $\DD$) if
$$
||\bmatrix q_{-1}&p_0\\
        p_0&q_0\endbmatrix||\le 1.
$$ 
Therefore, once the scattering matrix  
is not a  rational function of degree two,
then, for sure,  the holomorphic $S(t)$ does
not determine
$J $ uniquely.

The main result of this note is as follows.
\begin{theorem}\label{th3.1}
Given $S(t)$ from $H^\infty$, one can find an
inner function $\Phi(t)=\overline{\Phi(\bar
t)}$ such that 
\begin{equation*}
S_{\Phi}(t)=\bmatrix  s_-\Phi & s\\
         s& s_+ \overline \Phi
\endbmatrix(t),\qquad t\in \TT,
\end{equation*}
determines a Jacobi matrix of Szeg\"o class
uniquely. 
\end{theorem}
Let us point out that we did not
change the transmission coefficient $s$.
Moreover the new reflection coefficient
$s_-\Phi$ is still analytic one, so, the
corresponding Jacobi matrix $J$
coincides with $J_0$, when both are
restricted on
$l^2(\ZZ_{-}(-1))$.

Theorem \ref{th3.1} is the direct consequence
of the main result of \cite{Kh1} and the
following two lemmas.
\begin{lemma} Let the reflection coefficient
$s_-\in H^\infty$, and 
$s_-(\zeta)=a\zeta+\dots$, i.e.: $s_-(0)=0$.
Then \eqref{i16} holds if
\begin{equation}\label{i20}
\bmatrix
1&s_-\\
 0&s
\endbmatrix
\bmatrix
1&\bar t\\
\bar t&1
\endbmatrix=
\lim_{n\to\infty}\bmatrix
s & 0\\
s_+& 1
\endbmatrix
\bmatrix
f_n & \bar t g_n\\
\bar t f_n(\bar t)& g_n(\bar t)
\endbmatrix
 \end{equation}
with a suitable choice of sequences of
functions $f_n, g_n$ from
$H^2$, that is,  
$\frac{1+s_-/t} s\in H^2_{s_+}$ and 
$\frac{1+t s_-} s\in H^2_{s_+t^{-2}}$.
\end{lemma}
\begin{proof}
Since $s_-$ is holomorphic, we have
$\mathcal H_{s_-}=0$. Moreover, even
$\mathcal H_{s_-t^{-2}}\text{\bf 1}
=a\text{\bf 1}$, since $s_-(0)=0$. Therefore, 
\begin{equation}\label{i21}
K_{s_-}=\text{\bf
1}\quad \text{and}\quad
  K_{s_-t^{-2}}=\frac
{\text{\bf 1}}{\sqrt{1+a}}.
\end{equation}
We work with the second vector.
Because
of
\begin{equation}
\begin{split}
 \left\langle
 \bmatrix s&0\\s_+&1\endbmatrix
\bmatrix f^+(t)\\ \bar t f^+(\bar t)
\endbmatrix,
\bmatrix 1&s_-\\0&s\endbmatrix
\bmatrix  1\\ \bar t
\endbmatrix
\right\rangle =&
\left\langle
 \bmatrix s&0\\0&\bar s\endbmatrix
\bmatrix f^+(t)\\ \bar t f^+(\bar t)
\endbmatrix, \bmatrix 1\\ \bar t
\endbmatrix
\right\rangle\\
=&2(s f^+) (0),\quad f^+\in H^2,
\end{split}
\end{equation}
 and \eqref{i20}, the vector
$\frac{1+s_-/t} {s\sqrt{1+a}}$ 
 belongs to $H^2_{s_+}$
 and it is collinear with the reproducing
kernel of this space. Since its norm is one, 
we get
$$
K_{s_+}(\zeta)=\frac{1+s_-(\zeta)/\zeta}
{s(\zeta)\sqrt{1+a}}.
$$
Putting here $\zeta=0$, we get
$s(0)K_{s_+}(0)K_{s_-t^{-2}}(0)=1$. The proof
of the second identity is very similar.
\end{proof}

\begin{lemma}
Let the following approximation hold true
\begin{equation}\label{1}
\bmatrix
\bar t s_-\\
\frac{s-s(0)}t
\endbmatrix=
\lim_{n\to\infty}\bmatrix
s u_n\\
P_+ s_+ u_n
\endbmatrix,\qquad u_n\in H^2
\end{equation}
(we still assume that $s_-(0)=0$).
Then \eqref{i20} is also true.
 \end{lemma}

\begin{proof}
 We put
\begin{equation}\label{2}
\bar t f_n(\bar t)=
s(0)\bar t+ \bar t u_n(\bar t)
-P_- s_+ u_n
\end{equation}
and $g_n= f_n-(1-t^2)u_n$. In this case
\begin{equation*}
\bmatrix
f_n & \bar t g_n\\
\bar t f_n(\bar t)& g_n(\bar t)
\endbmatrix
=
\bmatrix f_n-u_n& t u_n\\
\bar t u_n(\bar t)& f(\bar t)-u(\bar t)
\endbmatrix
\bmatrix
1&\bar t\\
\bar t&1
\endbmatrix.
 \end{equation*}
In fact, we are going to show that
\begin{equation}\label{3}
\bmatrix
1&s_-\bar t\\
 0&s\bar t
\endbmatrix
 =
\lim_{n\to\infty}\bmatrix
s & 0\\
s_+& 1
\endbmatrix
\bmatrix
f_n-u_n &  u_n\\
\bar t u_n(\bar t)& \bar t \{f_n(\bar t)
-u_n(\bar t)\}
\endbmatrix.
 \end{equation}
Let us note that the equality in the
second column of \eqref{3} is just a
direct consequence of \eqref{1} and
\eqref{2}. Therefore, we need to check the 
 equality in the first column. In the first
entry of this column we have
\begin{equation*}
\begin{split}
\lim_{n\to\infty}
s(f_n-u_n)=&
\lim_{n\to\infty}
\{sf_n-s_-\bar t\}\\
=&
\lim_{n\to\infty}
\{s s(0)+ s u_n
-s P_+\bar t s_+(\bar t) u_n(\bar t)
-s_-\bar t\}\\
=&
\lim_{n\to\infty}
\{s s(0)+ s_-\bar t 
-s \overline{s_+}\bar t u_n(\bar t)
+s P_-\bar t s_+(\bar t) u_n(\bar t)
-s_-\bar t\}.
\end{split}
\end{equation*}
Using the unitary property of $S$--matrix
and \eqref{1}, we get
\begin{equation*}
\begin{split}
\lim_{n\to\infty}
s(f_n-u_n) 
=&
\lim_{n\to\infty}
\{s s(0) 
+s_- \overline{s}\bar t u_n(\bar t)
+s P_-\bar t s_+(\bar t) u_n(\bar t)
\}\\
=&
\lim_{n\to\infty}
\{s s(0) 
+s_- \overline{s_-} 
+s (\bar s-s(0))
\}\\
=&
\lim_{n\to\infty}
\{ 
|s_-|^2
+|s |^2
\}=1.
\end{split}
\end{equation*}
For the second entry we have
\begin{equation*}
\begin{split}
\lim_{n\to\infty}
\{s_+ (f_n -u_n) +\bar t u_n(\bar t)\} 
=&
\lim_{n\to\infty}
\{s_+ s(0)-s_+P_+\overline{s_+}\bar t
u_n(\bar t) + \bar t u_n(\bar t)\}\\
=&
\lim_{n\to\infty}
\{s_+ s(0)+s_+P_-\overline{s_+}\bar t
u_n(\bar t) + (1-|s_+|^2)\bar t u_n(\bar t)\}\\
=&
\lim_{n\to\infty}
\{s_+ s(0)+s_+ (\bar s-s(0)) + s\bar s\bar t
u_n(\bar t)\}\\ 
=&
s_+ s(0)+s_+ (\bar s-s(0)) + s
\overline{ s_-}=0.
\end{split}
\end{equation*}
\end{proof}

 \bibliographystyle{amsplain}

\end{document}